\documentclass{mcom-l}
\usepackage{amssymb}
\usepackage{amsmath}
\usepackage{color,xcolor}
\usepackage{graphicx}
\usepackage{bbm}

\numberwithin{equation}{section}

\begin{document}

\title{Newer sums of three cubes}

\author[S.G. Huisman]{Sander G. Huisman}
\address{Univ Lyon, Ens de Lyon, Univ Claude Bernard, CNRS, Laboratoire de Physique, F-69342 Lyon, France}
\email{s.g.huisman@gmail.com}
\thanks{The majority of the computations, totaling slightly above $10^5$ hours, were performed on a variety of machines of le p\^ole scientifique de mod\'elisation  num\'erique (PSMN) which is part of the \'Ecole Normale Sup\'erieure de Lyon. Stimulating and fruitful discussions with, and encouragement by T.~Browning, S.~Hervat, V.~Mathai, and E.~Trejo are acknowledged. The author also acknowledges the inspiration set forth by the YouTube video made by B. Haran on his channel Numberphile featuring T. Browning titled ``The Uncracked Problem with 33''. }


\date{\today}

\begin{abstract}
 The search of solutions of the Diophantine equation $x^3 + y^3 + z^3 = k$ for $k<1000$ has been extended with bounds of $|x|$, $|y|$ and $|z|$ up to $10^{15}$. The first solution for $k=74$ is reported. This only leaves $k=33$ and $k=42$ for $k<100$ for which no solution has yet been found. A total of 966 new solutions were found.
\end{abstract}

\maketitle

\section{Introduction}
The question of whether or not any $k\in \mathbbm{N}$ can be written as the sum of three cubes of integers:
\begin{align}
	x^3 + y^3 + z^3 = k \label{eq:stc}
\end{align}
has been investigated for centuries. Specific cases for $k$ can be excluded. For example, any $k \equiv \pm 4 \pmod{9}$ cannot have a solution since for any $n\in\mathbbm{N}$, $n^3 \equiv 0, \pm1 \pmod{9}$. In addition, for $k = 0$ the equation reduces to a specific case of Fermat's last theorem, which was proven by Euler to have no solutions besides the obvious trivial ones. It has not been proven that for all other $k$ solutions exist. Nor has there been any proof showing that any other $k$ has no solutions. Moreover, recent work of Colliot-Th\'el\`ene and Wittenberg \cite{colliot12} has shown that the integral Brauer-Manin obstruction is empty for these Diophantine equations. It is therefore conjectured that all other $k$ can be written as the sum of three cubes. 

For $k=1$ and $k=2$ (and therefore for $k=s^3$ and $k=2s^3$ for any $s\in\mathbbm{N}$) parametric solutions are known \cite{werebrusov08,mahler36,mordell42}:
\begin{align*}
	\left(9t^4\right)^3 + \left(3t-9t^4 \right)^3 + \left(1-9t^3 \right)^3 &= 1, \\
	\left(1 + 6t^3 \right)^3 + \left( 1 - 6t^3 \right)^3 + \left( -6t^2 \right)^3 &= 2.
\end{align*}
For this set of $k$ there are an infinite number of solutions that one can generate trivially. For other $k$ many solutions have been found for $k<1000$. The search for solutions accelerated in the '50s when digital computers became available. A complete historical overview and the progress over the years up to 2007 is given by \cite{beck07}. They end their paper with a list of $k$ for which no solution has been found yet: 3 values for $k<100$ and 27 values for $k<1000$. The most recent computer assisted search for solutions with bounds for $x$, $y$, and $z$ up to $10^{14}$ by Elsenhans and Jahnel from 2009 \cite{elsenhans09} still leaves 3 values for $k<100$ without solution and leave only 14 values for $k<1000$ without solution. In this paper the results for search bounds for $x$, $y$, and $z$ up to $10^{15}$ are presented.

\section{Methodology}
The numerical technique used is the same as the one used by Elsenhans and Jahnel \cite{elsenhans09}, and is based on the work of Elkies \cite{elkies00}. Here we look for large-number solutions of (\ref{eq:stc}). What that practically means is that there will be positive and negative numbers among $x$, $y$, and $z$. The problem can be rewritten as:
\begin{align}
	X^3 + Y^3 &= 1 + \left( \frac k z \right)^3, \label{eq:standardized}
\end{align}
where $X=x/z$ and $Y=y/z$. For large-number solutions for relatively small $k$ the last term is very small and effectively rational points are sought around the curve $Y^3 = 1 - X^3$. As the above equation is symmetric in $X$ and $Y$ the search space can be reduced by requiring $0<X<\sqrt[3]{1/2}$. Because the last term of (\ref{eq:standardized}) is nonzero solutions \textit{around} the line are sought for. The cubic curve is therefore covered by tiny `flagstones'\cite{elsenhans09} (called after the German word \textit{fliesen} meaning tiles) which are parallelograms that cover the curve, where 2 lines are parallel to the $Y$ axis and the other 2 sides follow the tangent of the curve. The width of these parallelograms depend on the search bound $B$ and the maximum size for $k$. Findings solutions of the form $x/z$, $y/z$ that fall inside these parallelograms is equivalent to looking for solution $x$, $y$, and $z$ inside a pyramid with sharp apex. In order to quickly iterate over lattice points inside this pyramid, lattice base reduction is applied to find basis vectors that are more accustomed to this geometry. For more details we refer to the documentation and the (German) comments inside the code of \cite{elsenhans09} and a more elaborate explanation of \cite{macleod11} for a very similar problem using identical technique.

\section{Results}
The results of \cite{elsenhans09} are taken as reference. With respect to this reference 966 new solutions have been found. The combined list now includes 15254 solutions. The number of solutions for each decade up to search bounds $B$ in the range from $10^2$ to $10^{14}$ have been roughly 1000 solutions per decade, which is in accordance with $\log(B)$\cite{heath92} scaling of the number of solutions up to bound $B$. The most exciting result is the discovery of the first solution for $k=74$:
\begin{align*}
	74 &= (-284650292555885)^3 + 66229832190556^3 + 283450105697727^3.
\end{align*}
For 3 values of $k$ a second solution has been found:
\begin{align*}
	606 &= (-170404832787569)^3 + (-16010802062873)^3 + 170451934224718^3, \\
	830 &= (-947922123009026)^3 + (-335912682279105)^3 + 961779444965911^3, \\
	966 &= (-1134209166959435)^3 + 291690681248788^3 + 1127741630138089^3.
\end{align*}
And for 9 values of $k$ for which only 2 solutions were known, now new solutions are available:
\begin{align*}
	30 &= (-662037799708799)^3 + 190809268841284^3 + 656711689254565^3, \\
	87 &= (-180751987188142)^3 + 40304620415495^3 + 180081502187630^3, \\
	327 &= (-170547806910404)^3 + (-55495503315494)^3 + 172484397062815^3, \\
	327 &= (-608821354548722)^3 + (-23497102374545)^3 + 608833020860980^3, \\
	402 &= (-577521772954300)^3 + 56488093736531^3 + 577341575648471^3, \\
	483 &= (-580121771294941)^3 + 247110530931002^3 + 564773663690516^3, \\
	735 &= (-857556020813401)^3 + 475896748252490^3 + 805621184663996^3, \\
	758 &= (-930395267860313)^3 + (-698249403713424)^3 + 1046417447600959^3, \\
	786 &= (-548172948160250)^3 + 377935686060853^3 + 480213828139669^3, \\
	831 &= (-210051134032889)^3 + (-160898634415892)^3 + 237716597381242^3.
\end{align*}
The question whether all $k\in \mathbbm{N}$ ($k \neq \pm4 \pmod{9}$, $k\neq0$) can be written as the sum of three cubes will likely remain an open question for many years to come. Our work does, however, give us new insight into the density of the solutions. The search bound was extended 10-fold since the last search \cite{elsenhans09}, leaving 13 unsolved $k<1000$ for which no solutions have been found yet: $33$, $42$, $114$, $165$, $390$, $579$, $627$, $633$, $732$, $795$, $906$, $921$, and $975$.

\bibliographystyle{alphanum}
\bibliography{paper}

\end{document}